\documentclass[12pt]{amsart} 

\begin{document}

\def\cal#1{{\mathcal #1}}
\def\R{{\mathbb R}} 
\def\tr{\mathop{\rm tr}}
\def\iint{\int\!\!\!\int}
\newcommand{\QED}{\hspace*{2em}\hfill$\bullet$}
\newtheorem{theorem}{Theorem}
\newtheorem{proposition}[theorem]{Proposition}
\newtheorem{definition}[theorem]{Definition}
\newtheorem{example}[theorem]{Example}
\newtheorem{lemma}[theorem]{Lemma}

\def\eqalign#1{\null\,\vcenter{\openup\jot \mathsurround=0pt
    \ialign{\strut\hfil$\displaystyle{##}$&$ \displaystyle{{}##}$\hfil
      \crcr#1\crcr}}\,}

\def\smallmatrix{\null\,\vcenter\bgroup \baselineskip=7pt
  \ialign\bgroup\hfil$\scriptstyle{##}$\hfil&&$\;$\hfil
  $\scriptstyle{##}$\hfil\crcr}
\def\endsmallmatrix{\crcr\egroup\egroup\,} \def\smatrix#1{\smallmatrix
  #1 \endsmallmatrix} \def\spmatrix#1{\big( \smallmatrix #1
  \endsmallmatrix \big)}

\title{Numerical integrators that contract volume}

\author{Robert I. McLachlan}
\address{\hskip-\parindent
        Robert I. McLachlan\\
        Mathematics\\
        Massey University\\
        Private Bag 11--222\\
        Palmerston North\\
		New Zealand}
\email{R.McLachlan@massey.ac.nz}

\author{G.R.W. Quispel}
\address{\hskip-\parindent
		G.R.W. Quispel\\
		Mathematics Department\\
		LaTrobe University\\
		Bundoora\\
		Melbourne 3083\\
		Australia}
\email{R.Quispel@latrobe.edu.au}

\thanks{Research at MSRI is supported in part by NSF grant DMS-9701755.}

\begin{abstract}
We study numerical integrators that contract phase space volume even
when the ODE does so at an arbitrarily small rate. This is
done by a splitting into two-dimensional contractive systems. 
We prove a sufficient condition for Runge-Kutta methods to have the
appropriate contraction property for these two-dimensional
systems; the midpoint rule is an example.
\end{abstract}

\maketitle

\subsection*{1. Introduction}

What is a dissipative system? In physics, the term usually refers
to possession of a scalar function (such as energy) which decreases 
in time, and one speaks of, e.g., the dissipative pendulum,
$\ddot x=-\sin x -\varepsilon\dot x$, for which $\frac{d}{dt}
(\frac{1}{2}\dot x^2
-\cos x) = -\varepsilon {\dot x}^2\le 0$. (See \cite{ma,mc-qu-ro1} for
some general formulations of such systems.) In dynamical systems, 
it usually refers to a decrease of phase space volume in time, as in the
``dissipative H\'enon map'' $(x,y)\mapsto (y, 1+bx-a y^2)$, with Jacobian
determinant $-b$---phase space area decreases if $|b|<1$. Another example
is the famous Lorenz system, which contracts volume at a constant rate.
In the numerical analysis of ODEs, it has been used to describe systems
that decrease some norm of the solution, either in the sense that
${d\over dt}\|x\|^2 < a-b\|x\|^2$ for some $a$, $b>0$,
or ${d\over dt}\|x\|^2 <0$ for all $\|x\|>R>0$ \cite{st-hu}.

In the field of geometric integration, much work has been done in 
maintaining the preservation of a conserved quantity (first integral)
\cite{gonzalez,mc-qu-ro1,mc-qu-ro2}, 
the decrease of a dissipated quantity (a Lyapunov function) 
\cite{mc-qu-ro1,mc-qu-ro2},
or the preservation of phase space volume \cite{fe-wa,mc-qu}. Here
we look at the missing case, and study how to maintain the property
of contracting phase space volume. 

Consider the ODE
\begin{equation}
\label{ode}
\dot x = f(x),\quad x\in \R^n
\end{equation}
with solution $x(t)$ and Jacobian (first variation) 
$A(t) = \partial x(t)/\partial x(0)$ which evolves according to
$$ \dot A = F A,\quad A(0)=I,$$
where $F(x)=df(x)$ is the derivative of the vector field $f$. We have 
$$\eqalign{{d\over dt} \det A  & = \det A \tr \left(A^{-1} \dot A\right) \cr
&=\det A \tr F}$$
so that phase space volume contracts, is preserved, or expands
when $\tr F<0$, $\tr F=0$, or $\tr F>0$ for all $x$, respectively. 
$\tr F$ is the divergence or trace of the vector field $f$.
Strongly contractive systems are those
for which there is a $b$ such that $\tr F<b<0$. In this case
any consistent numerical integrator will be contractive for small
enough time step $h$. Therefore we concentrate on weak contraction (defined
below), which
is a closed property and more difficult to preserve. 
It turns out that requiring contractivity for {\it all} $h>0$ 
and all contractive $f$ is prohibitively 
difficult, which leads to the following definition.
We consider one-%
step methods $x_{n+1} = g(x_n)$ with Jacobian $A = dg(x_0)$. 

\begin{definition}
The ODE (\ref{ode}) is (weakly) contractive
if $\tr F\le 0$ for all $x$. An integrator is (weakly) contractive if for
any matrix norm $\|\cdot\|$ and all
$L>0$ there is a time step $h^*>0$ such that $|\det A|\le 1$
for all $0<h<h^*$, for all $x$, 
and for all $f$ such that $\|F\|<L$ and $\tr F\le 0$.
\end{definition}

That is, there might be stiffness problems (for large $L$, $h^*$ might
be small), but the time step needed to preserve contractivity should not
tend to zero as $\tr F\to 0$. Note that a contractive integrator 
as defined here is not necessarily volume-preserving when the ODE is,
nor is the relative amount of contraction necessarily correct as 
$\tr F\to 0$. These would be true if we added the requirement
$\ln(\det A)/h\tr F \to 1$ uniformly as $\tr F\to0$ uniformly,
for all fixed $h<h^*$. The midpoint rule (see Proposition 3, below)
satisfies this, for example.

Since there are no known linearly covariant
volume-preserving schemes in more than
two dimensions \cite{fe-wa}, we expect that the same is true here,
and we immediately consider systems in two dimensions.

\subsection*{2. Dissipative schemes in two dimensions}

\begin{example} Euler's method is not contractive in two dimensions. 
\rm We have $x_{n+1} = x_n + h f(x_n)$ so $A= I + h F$. In two dimensions,
$$\det A = \det\left(
\begin{array}{cc}
1 + h F_{11} & h F_{12} \\
h F_{21} & 1 + h F_{22} \\
\end{array}
\right) 
= 1 + h \tr F + h^2 \det F.$$ 
So $\det A\le 1$ for all $h$ if $\det F \le 0$, and
$\det A\le 1$ for 
$$ h \le  {-\tr F \over \det F} $$
if $\det F>0$,
so small contractivity can require a small time step to be captured. 
\end{example}

Note that since $\det A = \det(I+h F)= \prod(1+h\lambda_i)$, where $\lambda_i$
are the eigenvalues of $F$, Euler's method {\it is} contractive in $n$
dimensions on systems with bounded negative eigenvalues. We
look at this further in Section 3.

\begin{proposition} The midpoint rule, 
$x_{n+1} = x_n + h f(\bar x)$, $\bar x = (x_n + x_{n+1})/2$, 
is contractive in two dimensions.
\end{proposition}
\begin{proof}
We have 
$$A = \Big(I-{1\over 2}h F(\bar x)\Big)^{-1} \Big(I+{1\over2}hF(\bar x)\Big),$$
so 
$$\det A = {1 + h e + h^2 d \over 1 - h e + h^2 d}$$
where $e = {1\over 2}\tr F(\bar x)$, $d = {1\over 4}\det F(\bar x)$.
Thus $(\det A)^2\le1$ if
$$ \left(1 + h e + h^2 d\right)^2 \le \left(1 - h e + h^2 d\right)^2 $$
or
$$ e \left(1 + h^2 d\right) \le 0 $$
Since $e\le 0$, this is true for all $h$ if $e=0$ (the well-known result
that the midpoint rule is area-preserving, or symplectic),
for all $h$ if $d\ge0$, or for $h<1/\sqrt{-d}$ if $d<0$.
\end{proof}

Proposition 3 can be generalized as follows.

\begin{proposition}
The symplectic Runge-Kutta methods with $b_i>0$ for all $i$
are contractive in two dimensions.
\end{proposition}
\begin{proof}
For the terminology, see \cite{sa-ca}. Our proof
closely follows their proof of symplecticity. An $s$-stage
 Runge-Kutta method is defined by
\begin{equation}
\label{rk1}
X_i = x_n+ h\sum_{j=1}^s a_{ij} f(X_j),
\end{equation}
\begin{equation}
\label{rk2}x_{n+1} = x_n + h \sum_{j=1}^s b_j f(X_j),
\end{equation}
and is symplectic if $b_i b_j - b_i a_{ij} - b_j a_{ji} = 0$ for all 
$i$ and $j$. Note that in two dimensions, $A^T J A = J \det A$, where
$J = \spmatrix{0 & 1 \cr -1 & 0}$, so we evaluate the left hand side. 
Let $D_i = df(X_i)= F(X_i)dX_i =: F_i A_i$. 
Differentiating (\ref{rk2}) gives
$$\eqalign{A^T J A & = \Big(I + h \sum_i b_i D_i\Big)^T J \Big(I+h \sum_j b_j D_j\Big) \cr
& = J + h \sum_i b_i \left(J D_i + D_i^T J\right) + h^2 \sum_{i,j} b_i b_j
D_i^T J D_j.}$$
Differentiating (\ref{rk1}) gives
\begin{equation}
\label{Ai}
 A_i = I + h\sum_j a_{ij} D_j
\end{equation}
or
$$ \eqalign{JD_i & = A_i^T J D_i - h \sum_j a_{ij} D_j^T J D_i}.$$
Inserting,
$$ A^T J A = J + h \sum_i b_i \left(A_i^T J D_i + D_i^T J A_i\right) 
+ h^2 \sum_{i,j} \left(b_i b_j - b_i a_{ij} - b_j a_{ji}\right) D_i^T J D_j.$$
The last term is zero because of the assumption on the coefficients
$b_i$, $a_{ij}$. Now $D_i = F_i A_i$, so
$$\eqalign{ A^T J A &= J + h\sum_i b_i A_i^T (J F_i + F_i^T J) A_i \cr
&= J + h \sum_i b_i A_i^T J A_i \tr F_i \cr
&= J \Big(1+ h \sum_i b_i \det A_i \tr F_i\Big) \cr
}$$
so
$$ 
\det A = 1 + h \sum_i b_i \det A_i \tr F_i.
$$
From (\ref{Ai}), $\det A_i$ is bounded and equal to $1+{\cal O}(h)$.
Using $b_i>0$ and $\tr F_i\le 0$ gives the result. 
\end{proof}

The assumption
$b_i>0$ is necessary. Suppose there are $s=2$ stages with $b_1>0$
and $b_2<0$. Then the vector $(b_i)$ lies in the fourth quadrant,
and all we know of the vector $(\tr F_i)$ is that it lies in the
third quadrant. In regions where the trace varies relatively quickly,
the angle between these two vectors can be less than ${\pi\over2}$,
leading to $(b_i)\cdot(\tr F_i)>0$ and $\det A>1$.

These methods actually preserve area when $\tr F=0$. In fact,
this is not necessary for contractivity in two dimensions, because
we can allow a small amount of ``numerical contractivity'' even
as $\tr F\to 0$; away from $\tr F=0$ the inherent contractivity
of the ODE contributes. It turns out that only methods
of order 2, 3, 6, 7,\dots, can achieve this.

\begin{lemma} Let $R(z)$ be the linear stability polynomial
of a consistent Runge-Kutta method. In two dimensions, the method
is contractive on linear ODEs if there is a $u^*>0$ such that
$$R(u)R(-u)\le 1,\quad R(i u)R(-i u)\le 1$$
for all $0\le u<u^*$.
\end{lemma}
\begin{proof}
In $n$ dimensions, a Runge-Kutta method on linear problems
$\dot x = Fx$ has derivative $A=R(hF)$. Therefore
$\det A = \prod_i R(h\lambda_i) = 1 + h\tr F + {\cal O}(h^2)$,
so the method is contractive if $\tr F<0$. If $\tr F=0$
we have to examine $\det A$ in more detail. In $n=2$
dimensions, there are only two such cases: the eigenvalues
can be $(u,-u)$ or $(i u,-i u)$. This gives the result.
\end{proof}

We note that the result also applies to nonlinear problems
with 1-stage methods, since then $F(x)$ is evaluated at
only a single point.

\begin{proposition}
\label{orderp}
Let the method have order $p$, so that
$R(z)=e^z + a z^{p+1}+b z^{p+2} + {\cal O}(z^{p+3})$.
If $4|(p+1)$ and $a<0$, or if $4|(p+2)$ and $b<a$, then the
method is contractive on linear problems in two dimensions.
\end{proposition}
\begin{proof}
We expand
$$ \eqalign{R(z)R(-z)-1 &= e^{-z}(a z^{p+1} + b z^{p+2})
+ e^z (a (-z)^{p+1} + b(-z)^{p+2}) +\dots \cr
& = a z^{p+1} (1-(-1)^p) + (b-a) z^{p+2}(1+(-1)^p) +\dots
}$$
The leading term must be negative for $z=u$ and for $z = iu$,
so it must be a fourth power. If $p$ is even, the leading
term is $z^{p+2}$ so $4|(p+2)$ and we need $b<a$; if
$p$ is odd, the leading term is $z^{p+1}$ so $4|(p+1)$
and we need $a<0$.
\end{proof}

An example is any 3-stage, 3rd order Runge-Kutta, which
has $R(z)=1+z+{1\over 2}z^2 + {1\over6}z^3$.

This result can be extended to more dimensions. For
example, a longer calculation shows that Proposition
\ref{orderp} holds with $p=3$ in three dimensions.
We are not sure how it extends to nonlinear systems.
It seems that if the eigenvalues of $F$ are varying rapidly,
contractivity could be lost. 

\subsection*{3. More than two dimensions}

For systems in more than two dimensions, we generalize the
volume-preserving method of Feng and Wang (\cite{fe-wa}; see
also \cite{mc-qu}). We write the ODE as a sum
of two-dimensional contractive systems 
(i.e., ones for which $\dot x_i=0$
except for two indices $i$), apply a contractive method
to each term, and compose the resulting maps with positive time steps. 
Since contractivity is
a semi-group property, we can build a contractive integrator of
order 1 or 2 in this way \cite{mc-qu}. 
This relies on the following proposition.

\begin{proposition}
Any $C^{r+1}$ contractive ODE is the sum of two-dimensional $C^r$
contractive ODEs.
\end{proposition}
\begin{proof}
Consider $\dot x = f(x)$, $F=df$. We shall
write $f$ in the form $f_i = \sum_j \partial_j L_{ij}$ (where $\partial_j
= \partial/\partial_{x_j}$.) 

Let $s_{ij}(x)$ be  $n^2$ functions with
$$ s_{ij}(x) + s_{ji}(x)\ge 0\quad{\rm and\ }
\sum_{i, j=1}^n s_{ij}(x)=1$$ 
for all $x$. Let 
$$ S_{ij} = \iint s_{ij}(x) \tr F(x)\, dx_i\, dx_j $$
where any values of the indefinite integrals can be taken. 
Let 
$$\tilde f_i = f_i - \sum_j\partial_j S_{ij} = f_i - \sum_j
\int s_{ij} \tr F\, dx_i,$$
so that
$$\tr d\tilde f = \Big(1-\sum_{i,j}s_{ij}\Big)\tr F = 0.$$
Thus, $\tilde f$ is traceless and can be written as
$$\tilde f_i = \sum_j \partial_j A_{ij}$$ 
where the matrix $A$ is antisymmetric and as smooth
as $f$ \cite{fe-wa,mc-qu}.
Therefore
\begin{equation}
\label{rep}
f_i = \sum_j \partial_j (A_{ij} + S_{ij})
\end{equation}
or $L=A+S$.

For an explicit splitting, we take the diagonal elements $s_{ii}(x)=0$.
Then $\dot x=f$ is the sum of the following $n(n-1)/2$ two-dimensional ODEs:
$$\eqalign{ \dot x_i &= \partial_j L_{ij} \cr
\dot x_j &= \partial_i L_{ji}\cr
\dot x_k &= 0 \hbox{\rm\ for\ }k\ne i,j\cr
}$$
for each pair $(i,j)$ of indices from $1$ to $n$.
Each is contractive because
each $A$ piece is traceless and each $S$ piece has trace
$(s_{ij} + s_{ji})\tr F\le0$.

One degree of smoothness is lost in this splitting,
because each $S$ piece depends on $\tr F$.
\end{proof}

An interesting solution is obtained by taking $s_{ii}=1/n$, $s_{ij}=0$
for $i\ne j$, and
$$nL_{ij} = \int f_i\, dx_j - \int f_j\, dx_i
+ \delta_{ij} \iint \tr F\, dx_i \, dx_i.$$
However, a more practical decomposition is to take the same $S$ but
$A_{ij}=0$ for $|i-j|>1$; this gives the minimum of $n-1$
two-dimensional ODEs. 

Although the above proof is constructive, 
it may be possible to find a more convenient splitting by ad hoc methods,
in some cases leading to an explicit contractive integrator.

Firstly,
if $f$ is the sum of {\it integrable} contractive vector fields, then
their flows can be composed to give a contractive integrator for $f$.
For example, the Lorenz system,
$$
\dot x = \left(\begin{array}{ccc}
-\sigma & \sigma & 0 \\
\rho & -1 & 0 \\
0 & 0 & -\beta
\end{array}
\right)x
+ \left( \begin{array}{ccc}
0 & 0 & 0 \\
0 & 0 & -x_1 \\
0 & x_1 & 0 \\
\end{array}\right)
\nabla\left( {x_2^2 + x_3^2\over 2}\right),
$$
is the sum of a linear, contractive part and a Poisson, non-contractive
part, each of which may be solved exactly, giving an integrator with
{\it exactly} correct contractivity.

Secondly, it may be possible to use a simpler method, such
as Euler, on some of the pieces. Here are some criteria
which allow this.

\begin{proposition}\label{euler_n}
Euler's method is contractive in $n$
dimensions if there is a $b$ such that $\tr(F^2)>b>0$. This condition is
equivalent to $\|S\|^2 > \|A\|^2+b$, where
$F=A+S$, $A=-A^T$, $S=S^T$, and $\|\cdot\|$ is the
Frobenius (sum of squares) norm. This condition is satisfied
if all the eigenvalues of $F$ are bounded away from the sectors
${\pi\over4}<|\theta|<{3\pi\over4}$; in particular,
if they are all real and bounded away from zero.
\end{proposition}
\begin{proof}
Let $\lambda_i$ be the eigenvalues of $F$. For Euler's
method we have
$$\eqalign{
\ln \det A &= \ln \det (I + h F) \cr
&= \ln \prod_i (1 + h \lambda_i) \cr
&= \sum \ln (1 + h\lambda_i) \cr
&= h\sum_i \lambda_i - {1\over 2} h^2 \sum_i \lambda_i^2
+{\cal O}(h^3) \cr
&=h\tr F-{1\over2}h^2\tr(F^2)+
{\cal O}(h^3). \cr
}$$
If there is a $b$ such that $\tr(F^2)>b>0$, this is less than $0$ for
all small enough $h$, i.e., the method is contractive.
Splitting $F$ into
its symmetric and antisymmetric parts, 
$$ \tr(F^2) = \sum_{i,j} F_{ij} F_{ji} 
= \sum_{i,j} (S_{ij} + A_{ij})(S_{ij} - A_{ij})
= \|S\|^2 - \|A\|^2,$$
giving the second part of the proposition.
Now $\tr(F^2)=\sum_i \lambda_i^2$, and if each $\lambda_i$ is
outside the specified sectors, then each real eigenvalue
or complex conjugate pair of eigenvalues gives a positive
contribution to this sum, giving the last part of the proposition.
\end{proof}

Note that the eigenvalues of elliptic or nearly elliptic
fixed points lie near the imaginary axis---right in the middle 
of the ``bad'' sector. Perhaps this was only to be expected.

Experts will recognize the last part of Proposition \ref{euler_n}
as the appearance of an {\it order star} of a Runge-Kutta
method \cite{is-no} (the set $\{z:|R(z)|<|e^z|$ where 
$R(z)$ is the method's
linear stability polynomial). For linear problems,
or nonlinear problems with 1-stage methods, a method is 
more contractive than the flow of the ODE if $h\lambda_i$
lies in the order star of the method for each eigenvalue
$\lambda_i$. However, this seems rather restrictive so
we do not explore further.

\begin{proposition}
There are explicit contractive integrators.
\end{proposition}
\begin{proof}
Let $f$ be any contractive vector field with $\|F\|<L$. Because
eigenvalues vary continuously and can only become imaginary when
two eigenvalues meet, and because symmetric matrices
have real eigenvalues,
there is a symmetric, traceless matrix $M$ with distinct
eigenvalues such that the derivative of $f_1:=f-Mx$
has real eigenvalues. Let $f_2 := Mx$ and split $f=f_1+f_2$.
$f_1$ is contractive and admits an explicit contractive
integrator (e.g. Euler's method, see Proposition 8); 
$f_2$ is traceless and can be solved explicitly. Composing
these maps gives the result.
\end{proof}

We close with some open questions we hope to report on in the future.
\begin{enumerate}
\item Are there explicit contractive integrators of any order?
(Proposition 9 constructs a first order method.) There are
if one only demands {\it linear} contractivity. The order 
cannot be increased by composition, because the adjoint of
Euler's method---backward Euler---is not contractive for $f_1$.
\item The present method reduces to
the volume-preserving method of Feng and Wang \cite{fe-wa}
when the vector field $f$ is traceless.
There is another approach to volume-preserving integration
due to Quispel \cite{quispel} and to 
Shang \cite{shang}, which does not rely on a splitting at all;
moreover, it has a generalization to systems preserving
non-Euclidean measures, which we have not even considered
here. Can this approach be carried over to the contractive case?
\item The splitting used in the proof of Proposition 5
writes $f=a+b$ where $\tr(da)=0$  and $b\equiv 0$ when $\tr(df)\equiv 0$.
Are there splittings with the property that $b(x)=0$ when $\tr(df(x))=0$?
If so, they could be used for systems in which $\tr(df(x))$ 
changes sign on a compact hypersurface; the interior would 
then be invariant and one could construct an integrator which 
preserved it and was contractive there. This was done for
the case of dissipation of scalar functions in \cite{mc-qu-ro2}.
\end{enumerate}

\subsection*{Acknowledgements}

The authors thank John Butcher and all the organizers of ANODE for 
financial support and for providing the atmosphere in which this 
paper was conceived, and the MSRI where it was concluded. We
also thank the Marsden Fund of the Royal Society of New
Zealand for their financial support.

\end{document}